\documentclass[11pt,a4paper]{amsart}
\usepackage[T1]{fontenc}
\usepackage[utf8]{inputenc}
\usepackage{lmodern,microtype}
\usepackage[margin=29mm]{geometry}
\usepackage{amsmath,amssymb,amsthm,mathtools}
\usepackage{enumitem,booktabs}
\usepackage[hidelinks]{hyperref}
\usepackage{xurl}
\hypersetup{
  pdftitle={On homogeneous hypertournaments},
  pdfauthor={Matěj Konečný}
}
\newtheorem{theorem}{Theorem}[section]
\newtheorem{lemma}[theorem]{Lemma}
\newtheorem{corollary}[theorem]{Corollary}
\newtheorem{proposition}[theorem]{Proposition}
\newtheorem{observation}[theorem]{Observation}
\newtheorem{conjecture}[theorem]{Conjecture}
\newtheorem{question}[theorem]{Question}
\theoremstyle{definition}
\newtheorem{definition}[theorem]{Definition}
\theoremstyle{remark}
\newtheorem{remark}[theorem]{Remark}

\newcommand{\K}{\mathcal K}
\newcommand{\F}{\mathcal F}
\newcommand{\Forb}{\operatorname{Forb}}
\newcommand{\parity}{\operatorname{par}}

\newcommand{\str}{\mathbf}
\title[On homogeneous hypertournaments]{On homogeneous hypertournaments}
\author[M. Kone\v cn\'y]{Mat\v ej Kone\v cn\'y}
\address{Charles University, Faculty of Mathematics and Physics\\Department of Applied Mathematics (KAM)\\
Prague, Czech Republic}
\email{matej@kam.mff.cuni.cz}
\date{}

\begin{document}
\begin{abstract}
Sullivan recently constructed continuum many countable homogeneous
$n$-hypertournaments for every $n\ge3$ by forbidding configurations
he calls pair-cancelling. We show that, for every family $\mathcal F$
of pair-cancelling $n$-hypertournaments, the class of all linearly
ordered finite $\mathcal F$-free $n$-hypertournaments is Ramsey and
has the ordering property. In particular, freely ordered $\str{H}_4$-free
$3$-hypertournaments are Ramsey, answering a question of Cherlin,
Hubi\v cka, Ne\v set\v ril, and the author from 2021.
For every odd $n\ge3$, we give a simple Ramsey expansion with
the expansion property for even $n$-hypertournaments by naming
an $(n-1)$-hypertournament and a free linear order.

We also show that, for $n\ge6$, no unbounded hereditary strong amalgamation class of finite $n$-hypertournaments has EPPA. Moreover, we prove the failure of EPPA for all Sullivan's classes with $n\ge4$, and for the even $5$-hypertournaments.

We conjecture that even
hypertournaments and classes omitting a family of pair-cancelling hypertournaments, together with the homogeneous cyclic order in
arity three, are the only countably infinite homogeneous
$n$-hypertournaments for $n\ge3$. We give some supporting evidence.
\end{abstract}

\maketitle

\section{Introduction}\label{sec:intro}
A structure is \emph{homogeneous} if every
isomorphism between finite induced substructures extends to an
automorphism. An $n$-hypertournament assigns one of the two alternating orientations
to each $n$-element set. A 2-hypertournament is simply a tournament.

The countably infinite homogeneous tournaments were
classified by Lachlan~\cite{lachlan1984countable}; there are three up to
isomorphism. The four $4$-constrained homogeneous 3-hypertournaments were described
by Cherlin, Hubička, Nešetřil, and the author in our EuroComb
abstract~\cite{cherlin2021ramsey}; Cherlin also classified finite
homogeneous hypertournaments~\cite{Cherlin2019Notes}.
Cherlin~\cite[Appendix~B]{Cherlin2022} asked for a classification of homogeneous
$3$-hypertournaments (Problem~23) and whether infinitely many homogeneous
$n$-hypertournaments exist (Problem~24).

Sullivan~\cite{Sullivan2026} recently answered the abundance question:
there are $2^{\aleph_0}$ nonisomorphic countable homogeneous
$n$-hypertournaments for every $n\ge3$. He constructs them by
forbidding finite configurations he calls \emph{pair-cancelling}
(Definition~\ref{def:pair-cancelling}).

In~\cite{cherlin2021ramsey}, we
studied Ramsey expansions of 4-constrained 3-hypertournaments and
left the $\str{H}_4$-free ones as the only open case. This was further discussed in~\cite{hubicka2025twenty,Konecny2023phd}.
Sullivan also asks whether this
structure and his new examples have precompact Ramsey expansions.
We answer these questions with the following theorem.

\begin{theorem}
\label{thm:intro-link}
Let $n\ge3$, let $\F$ be a family of finite
pair-cancelling $n$-hypertournaments, and put $\K=\Forb(\F)$.
The expansion $\K^<$ of $\K$ by free linear orders is Ramsey
and has the ordering property.
\end{theorem}

In particular, the $\str{H}_4$-free case is included, since $\str{H}_4$ is pair-cancelling.
In~\cite{cherlin2021ramsey,hubicka2025twenty} we incorrectly inferred from non-completable
partial structures that the Hubi\v cka--Ne\v set\v ril theorem~\cite{Hubicka2016} could
not apply directly to the freely ordered class. We explain the
mistake in Remark~\ref{rem:bicycles}.

For every odd $n\ge3$, there is another natural class: the class of
\emph{even} $n$-hypertournaments. Its Fra\"iss\'e limit is obtained from
the generic $(n-1)$-hypertournament by a parity condition on $n$-subsets.
These classes were first identified by Cherlin~\cite{Cherlin2019Notes}, and the ternary case
was discussed in~\cite{cherlin2021ramsey}, where we found its Ramsey expansion.
Here we extend the construction and prove (see Section~\ref{sec:potential} for definitions):
\begin{theorem}\label{thm:potential}
For every odd $n\ge3$, the class $\mathcal E_n^*$
is Ramsey and has the expansion property relative
to $\mathcal E_n$.
\end{theorem}

\medskip

We also study the extension property for partial automorphisms (EPPA) and obtain the following negative results (see Section~\ref{sec:eppa} for definitions and more results):

\begin{theorem}\label{thm:intro-eppa}
\leavevmode
\begin{enumerate}[label=(\roman*),nosep]
\item For $n\ge6$, no unbounded hereditary strong amalgamation
class of finite $n$-hypertournaments has EPPA.
\item For $n\ge4$, if $\F$ is any family of finite pair-cancelling
$n$-hypertournaments, then $\Forb(\F)$ does not have EPPA.
\item The class of finite even $5$-hypertournaments does not have EPPA.
\end{enumerate}
\end{theorem}

\medskip

Finally, in Section~\ref{sec:classification-evidence}, we discuss classification of countably infinite homogeneous hypertournaments. We propose the following conjecture:
\begin{conjecture}\label{conj:all-homogeneous}
Let $n\ge3$ and let $\str{M}$ be a countably infinite homogeneous
$n$-hypertournament. Then exactly one of the following holds:
\begin{enumerate}[label=(\roman*),leftmargin=*,itemsep=1pt]
\item $\operatorname{Age}(\str{M})=\Forb(\F)$ for a family $\F$ of finite
pair-cancelling $n$-hypertournaments;
\item $n$ is odd and $\operatorname{Age}(\str{M})$ is the class of even
$n$-hypertournaments;
\item $n=3$ and $\str{M}$ is the countable homogeneous cyclic order.
\end{enumerate}
\end{conjecture}

This is supported by complete computer-assisted classifications of the unbounded $(n+1)$-constrained
$n$-hypertournament classes, for all $n$, and of the $3$-hypertournament classes described by forbidding structures of size at most $5$, see Section~\ref{sec:classification-evidence}.

If the conjecture holds,
every countably infinite homogeneous $n$-hypertournament
($n\ge3$) has a precompact Ramsey expansion in a finite relational
language: use Theorems~\ref{thm:intro-link} and~\ref{thm:potential} for the first two cases, and a compatible linear order for
the cyclic order~\cite{cherlin2021ramsey}. In sharp contrast, no countably infinite homogeneous
hypertournaments are presently known to have EPPA. 

\medskip

\paragraph{\textbf{AI disclaimer.}}
This work was developed with iterative use of ChatGPT (GPT-5.6 Sol and GPT-6),
both in the mathematical research and in drafting and editing the manuscript.
I have directed, reviewed and verified the work in full and take responsibility
for all its contents.

\medskip

\paragraph{\textbf{Formal verification.}}
All results of this paper have been formalised in Lean~4; see \url{https://github.com/shulda/hypertournaments}.
A \href{https://github.com/shulda/hypertournaments/wiki/Main-article}{statement-by-statement correspondence}
and a \href{https://github.com/shulda/hypertournaments/blob/main/docs/theorem-map.json}{machine-readable theorem map}
are available there. The \href{https://github.com/shulda/hypertournaments/tree/97ece7669d91ffeff83227cc2d2b8cd64f141bfd}{Lean source snapshot}
fixes the formalisation version accompanying this article.
The Lean formalisation was developed almost entirely
autonomously by AI, primarily as an additional check of the proofs, without
aiming to meet the usual standards of style and reusability in formalised
mathematics.

\section{Preliminaries}\label{sec:prelim}
\subsection{Hypertournaments}
Generalising tournaments, Cherlin introduced $n$-hyper\-tour\-na\-ments as structures with one injective $n$-ary relation such that the automorphism group of every subset of $n$ vertices is the alternating group $\mathrm{A}_n$~\cite{Cherlin2019Notes,Cherlin2022}. In this note, it will be convenient to treat $n$-hyper\-tour\-na\-ments as structures with an alternating $\mathbb Z_2$-function on injective $n$-tuples as defined below:

We write $\mathbb Z_2=\{0,1\}$, with addition taken modulo~$2$. Denote by
$\parity(\pi)\in\mathbb Z_2$ the parity of a permutation $\pi$.
Throughout the paper we will denote structures by bold letters such as $\str X$ and their vertex sets by normal letters ($X$).
For a set $X$ and $k\ge1$, write $X^{\underline{k}}$
for the set of injective ordered $k$-tuples from $X$.
A function $f\colon X^{\underline{k}}\to\mathbb Z_2$
is \emph{symmetric} or \emph{alternating}, respectively,
if for every permutation $\pi$ and every injective
$k$-tuple $(x_1,\ldots,x_k)$,
\[
 \begin{aligned}
 f(x_{\pi(1)},\ldots,x_{\pi(k)})
   &=f(x_1,\ldots,x_k)
       &&\text{(symmetric)},\\
 f(x_{\pi(1)},\ldots,x_{\pi(k)})
   &=f(x_1,\ldots,x_k)+\parity(\pi)
       &&\text{(alternating)}.
 \end{aligned}
\]
The latter is the additive encoding of ordinary
sign alternation: after identifying $a\in\mathbb Z_2$
with $(-1)^a\in\{+1,-1\}$, odd permutations change
the sign. In particular, a symmetric function
is equivalently a function on $\binom Xk$.
Note that the sum
of two alternating functions on the same domain is symmetric, while
adding a symmetric function to an alternating
one gives another alternating function.

For $n\ge2$, an \emph{$n$-hypertournament} $\str{X}$ consists of a
vertex set $X$ (finite or infinite) together with an alternating
\emph{orientation map}
\[
 \omega_{\str{X}}\colon X^{\underline n}
                       \longrightarrow\mathbb Z_2.
\]
Equivalently, $\str{X}$ can be viewed as a relational structure with one
$n$-ary relation
\[
 R^{\str{X}}=\{(x_1,\ldots,x_n):
                         \omega_{\str{X}}(x_1,\ldots,x_n)=0\}.
\]
We switch freely between the relational and additive presentations.
Unless explicitly stated otherwise, classes of hypertournaments in this
paper consist of finite structures.
For a class $\K$
of $n$-hypertournaments, $\K^<$ denotes its expansion by \emph{free
linear orders}, which means that every linear order is permitted.

An \emph{embedding} is an injection preserving every orientation.
For a family $\F$ of finite $n$-hypertournaments, $\Forb(\F)$ denotes
the class of all finite $n$-hypertournaments with no embedding from any member of $\F$.

\subsubsection{The $(n+1)$-vertex types}\label{sec:local-types}
Let $\str X$ be an $n$-hypertournament on $n+1$ vertices, and enumerate
$X = x_0 < \cdots < x_n$ (the order is not part of the structure).
For every $0\leq i\leq n$, let $s_i\in\mathbb Z_2$ be the orientation of the
$n$-tuple omitting $x_i$, and put $t_i=s_i+i$.
Put
\[
 d(\str X)=\min\bigl(|\{i:t_i=0\}|,|\{i:t_i=1\}|\bigr).
\]
Note that $0\le d(\str X)\le q=\left\lfloor\frac{n+1}{2}\right\rfloor$. We will show below that the isomorphism type of $\str X$ only depends on $d(\str X)$;
following Cherlin's classification of these small types
\cite{Cherlin2019Notes}, we denote this type by
$\str{T}_{n,d}$. For $n=2$, the two types
$\str{T}_{2,0}$ and $\str{T}_{2,1}$ are, respectively,
the directed $3$-cycle and the transitive $3$-vertex
tournament. When $n=3$, the
three types $\str{T}_{3,0},\str{T}_{3,1},
\str{T}_{3,2}$ are usually called $\str{H}_4,\str{O}_4,\str{C}_4$.
For $\str{H}_4$, all $t_i$ are equal; for $\str{O}_4$,
exactly one $t_i$ differs from the other three; and for
$\str{C}_4$, each of the two values occurs twice.

If $\F$ consists of $(n+1)$-vertex $n$-hypertournaments,
we say that $\Forb(\F)$ is \emph{$(n+1)$-constrained}.

\begin{lemma}\label{lem:local-type-isomorphism}
Let $\str X$ and $\str X'$ be $n$-hypertournaments with $n+1$ vertices, for some $n\geq 2$. Then $d(\str X) = d(\str X')$ if and only if $\str X$ and $\str X'$ are isomorphic.
\end{lemma}
\begin{proof}
Fix an enumeration $(x_0,\ldots,x_n)$ of $X$,
and let $t_0,\ldots,t_n$ be defined as above.
For every permutation $\pi$ of $\{0,\ldots,n\}$,
let $t_i^\pi$ be the value defined in the same way
from the enumeration $(x_{\pi(0)},\ldots,x_{\pi(n)})$.
Then
\[
 t_i^\pi=t_{\pi(i)}+\parity(\pi)
 \qquad(0\le i\le n).
\]
Indeed, after omitting $x_{\pi(i)}$, the induced
permutation of the remaining $n$ vertices has parity
$\parity(\pi)+i+\pi(i)$. Thus the corresponding
orientation is $s_{\pi(i)}+\parity(\pi)+i+\pi(i)$;
adding $i$ gives the displayed formula. Consequently,
$d(\str X)$ is an isomorphism invariant.

Conversely, suppose that $d(\str X)=d(\str X')$.
By identifying their vertex sets, we may assume that $X=X'$
and compute the vectors $(t_0,\ldots,t_n)$ and
$(t'_0,\ldots,t'_n)$ using the same enumeration
$(x_0,\ldots,x_n)$.

As $d(\str X) = d(\str X')$, there is a permutation $\pi$ of $\{0,\ldots,n\}$
and $b\in\mathbb Z_2$ such that $ t'_i=t_{\pi(i)}+b$, for every $0\leq i\leq n$.
If $\parity(\pi)=b$, then the map
$x_{\pi(i)}\mapsto x_i$ is an isomorphism
$\str X\to\str X'$.
Otherwise choose distinct indices $r,s$ with $t_r=t_s$,
which exist since $n+1\ge3$. Replacing $\pi$ by
$(r\,s)\circ\pi$ preserves all equalities
$t'_i=t_{\pi(i)}+b$ but reverses the parity of $\pi$,
giving the desired isomorphism.
\end{proof}

\subsection{Amalgamation}\label{sec:prelim-amalgamation}
For a structure $\str{M}$, its \emph{age} $\operatorname{Age}(\str{M})$ is the class
of all finite structures which embed into $\str{M}$.
Let $\K$ be a class of finite structures in a fixed relational language.
The class $\K$ is \emph{hereditary} if it is closed under induced
substructures, and has the \emph{joint embedding property} 
if any two members embed into a third. Given embeddings
$\alpha_i\colon \str{A}\to \str{B}_i$, $i=1,2$, an \emph{amalgam} is a structure
$\str{C}$ together with embeddings $\beta_i\colon \str{B}_i\to \str{C}$ satisfying
$\beta_1\alpha_1=\beta_2\alpha_2$. The \emph{amalgamation property} asserts the existence of such a $\str{C}$ for every diagram. The
amalgam is \emph{strong} if
\[
 \beta_1[\str{B}_1]\cap\beta_2[\str{B}_2]
   =\beta_1\alpha_1[\str{A}].
\]
If all diagrams admit strong amalgams, then $\K$ has the \emph{strong
amalgamation property}. A hereditary class with the joint embedding property and the amalgamation property
is called an \emph{amalgamation class}.

Fra\"iss\'e's theorem~\cite{Fraisse1953,Fraisse1986} asserts that
an amalgamation class in a finite relational language is the age of a
unique countable homogeneous structure, called its
\emph{Fra\"iss\'e limit}. The limit is infinite precisely when the class
has members of arbitrarily large finite size (in which case we call
the class \emph{unbounded}).

\subsection{Ramsey classes}\label{sec:prelim-ramsey}
For structures $\str{A}$ and $\str{B}$, let $\operatorname{Emb}(\str{A},\str{B})$ denote the
set of all embeddings $\str{A}\to \str{B}$. Given $\str{A},\str{B},\str{C}$ in the same language and
an integer $r\ge2$, the notation
\[
 \str{C}\longrightarrow(\str{B})^{\str{A}}_r
\]
means that for every colouring
$\chi\colon\operatorname{Emb}(\str{A},\str{C})\to\{1,\ldots,r\}$ there is an
embedding $f\in\operatorname{Emb}(\str{B},\str{C})$ such that $\chi(f\circ e)$ is
constant for $e\in\operatorname{Emb}(\str{A},\str{B})$. A class $\K$ is
\emph{Ramsey} if for all $\str{A},\str{B}\in\K$ and $r\ge2$ there is $\str{C}\in\K$
with $\str{C}\longrightarrow(\str{B})^{\str{A}}_r$.
The classical Ne\v set\v ril--R\"odl theorem
\cite{Nevsetvril1977,nevsetvril1977partitions} says that the class
of all finite linearly ordered structures in a fixed finite relational
language is Ramsey (this was also proved independently by Abramson and Harrington~\cite{Abramson1978}).

Let $L\subseteq L^+$ be relational languages. The \emph{$L$-reduct}
of an $L^+$-structure is obtained by forgetting the relations in
$L^+\setminus L$. An \emph{expansion} $\K^+$ of a class $\K$ of
$L$-structures is a class of $L^+$-structures whose $L$-reducts are
precisely the members of $\K$; it is a \emph{Ramsey expansion} if
$\K^+$ is Ramsey.

Following Nguyen Van Th\'e~\cite{The2013}, we say that the
expansion has the \emph{expansion property} relative to $\K$ if
for every $\str{A}\in\K$ there exists $\str{B}\in\K$ such that every expansion
of $\str{A}$ embeds into every expansion of $\str{B}$. For expansions by free
linear orders, this is the \emph{ordering property}, already studied
by Ne\v set\v ril and R\"odl~\cite{nevsetril1975type,Nesetvril1978}.

The connection between Ramsey classes and topological dynamics
discovered by Kechris, Pestov, and Todor\v cevi\'c~\cite{Kechris2005} in 2005
revitalised structural Ramsey theory and led to a systematic search
for Ramsey expansions. Nowadays, Ramsey expansions of many classes
of finite structures are well understood
(see e.g.~\cite{Nevsetvril1976,Nevsetvril2007,Bodirsky2010,laflamme2010partition,
sokic2013ramsey,Jasinski2013,sokic2015semilattices,Aranda2017a,
Hubicka2017designs,bhat2016ramsey,Hubivcka2014,Hubicka2016,Konecny2018b,
Aranda2017c,Hubicka2017sauerconnant,Aranda2017,HubickaSauerSmetric,
Lee2026CRelations,Sam});
see~\cite{hubicka2025twenty} for a survey.
Another important development in structural Ramsey theory is the recent
girth Ramsey theorem~\cite{reiher2023girth}.

Most of the currently known Ramsey classes can be obtained using
the Hubi\v cka--Ne\v set\v ril theorem~\cite{Hubicka2016}.
A relational structure is \emph{irreducible} if every two
distinct vertices occur together in a relation tuple. A
\emph{homomorphism-embedding} $f\colon \str{C}\to \str{D}$ is a homomorphism whose
restriction to every irreducible substructure of $\str{C}$ is an embedding.
A \emph{completion} of $\str{C}$ in $\K$ is a homomorphism-embedding
$\str{C}\to \str{C}'$ for some $\str{C}'\in\K$; it is \emph{strong} if the map is
injective.

We recall the following definition from
\cite[Definition~2.8]{Hubicka2016}.
\begin{definition}
\label{def:locally-finite}
Let $\mathcal R$ be a class of finite irreducible structures,
and let $\K\subseteq\mathcal R$.
We call $\K$ a \emph{locally finite subclass} of $\mathcal R$ if,
for every $\str{C}_0\in\mathcal R$ there exists an integer
$n=n(\str{C}_0)\ge1$ such that every finite structure $\str{C}$
has a completion in $\K$ whenever
\begin{enumerate}[label=(LF\arabic*),nosep]
\item\label{item:LF1} $\str{C}_0$ is an $\mathcal R$-completion
of $\str{C}$;
\item\label{item:LF2} every irreducible substructure of $\str{C}$
belongs to $\K$;
\item\label{item:LF3} every substructure of $\str{C}$ on at most
$n(\str{C}_0)$ vertices has a completion in $\K$.
\end{enumerate}
\end{definition}

We can now state the Hubi\v cka--Ne\v set\v ril theorem~\cite[Theorem~2.11]{Hubicka2016}:
\begin{theorem}
\label{thm:HN}
Let $\mathcal R$ be a Ramsey class of finite irreducible
relational structures. If $\K\subseteq\mathcal R$ is hereditary,
locally finite, and has the strong amalgamation property,
then $\K$ is Ramsey.
\end{theorem}

\section{Omitting pair-cancelling hypertournaments}\label{sec:links}
In this section we will prove Theorem~\ref{thm:intro-link}. We will assume $n\geq 3$ throughout.

Let $\str X=(X,\omega_\str X)$ be an $n$-hypertournament and let $u,v\in X$ be distinct vertices of $\str X$. The \emph{link} of $u$ and $v$ is the function from $(X\setminus\{u,v\})^{\underline{n-2}}$ to $\mathbb Z_2$ given by
\[
 (x_1,\ldots,x_{n-2})\longmapsto
 \omega_{\str{X}}(u,v,x_1,\ldots,x_{n-2}).
\]
An odd permutation of the tuple adds $1$ to its link value, and reversing
$u,v$ adds $1$ to every link value.

The following definition is due to Sullivan~\cite{Sullivan2026}:
\begin{definition}
\label{def:pair-cancelling}
An $n$-hypertournament $\str{X} = (X,\omega_\str X)$ with at least $n+1$ vertices is
\emph{pair-cancelling} if, for every pair of distinct vertices $u,v$
and every linear order $\prec$ on $X\setminus\{u,v\}$,
the values $\omega_{\str{X}}(u,v,x_1,\ldots,x_{n-2})$
on $\prec$-increasing $(n-2)$-tuples $(x_1,\ldots,x_{n-2})$ are not all equal.
\end{definition}
For $n=3$, this says that for each $u,v$ both values $0$ and $1$
occur among $\omega_{\str{X}}(u,v,w)$ as $w$ varies.
Equivalently, a hypertournament is pair-cancelling if no pair has
a link whose increasing tuples all have the same value in some order.

\begin{lemma}\label{lem:local-types}
For $n\ge3$, an $(n+1)$-vertex
$n$-hypertournament is pair-cancelling if and
only if its type is among
$\str{T}_{n,0},\ldots,\str{T}_{n,q-2}$, where
$q=\lfloor(n+1)/2\rfloor$.
\end{lemma}
\begin{proof}
Let $\str{X}$ be an $(n+1)$-vertex
$n$-hypertournament. Fix distinct vertices $u,v$
and list the remaining vertices as $y_0,\ldots,y_{n-2}$. Use the order
$u<v<y_0<\cdots<y_{n-2}$ to compute the values $t_i$ from Section~\ref{sec:local-types}.
Writing $\widehat{y_j}$ for omission of $y_j$, we have
\[
 \omega_{\str{X}}(u,v,y_0,\ldots,\widehat{y_j},\ldots,y_{n-2})
   =t_{j+2}+j\qquad(0\le j\le n-2).
\]
Thus the increasing tuples in the link all have the
same value exactly when $t_2,\ldots,t_n$ alternate
between $0$ and $1$. Reordering the $y_j$ permutes
these $n-1$ entries, possibly adding $1$ to all of them.
Consequently, some order makes the link constant exactly when
the numbers of zeroes and ones among the remaining
$n-1$ entries differ by at most one. Such entries can be
obtained by deleting two entries from $(t_0,\ldots,t_n)$ exactly
when $d(\str X)\ge q-1$. Hence $\str{X}$ is pair-cancelling exactly when
$d(\str X)\le q-2$.
\end{proof}

\subsection{Completions and amalgamation}
We first show how to complete structures in the language $L=\{R,<\}$.

\begin{proposition}\label{prop:link-completion}
Let $n\ge3$, let $\F$ be a family of finite $n$-hypertournaments,
each of which is pair-cancelling, and put
$\K=\Forb(\F)$. Suppose that a finite $L$-structure $\str{C}$
satisfies:
\begin{enumerate}[label=(\roman*),ref=(\roman*),nosep]
\item\label{cond:completion-order} the relation $<^{\str{C}}$ admits a linear extension, and
\item\label{cond:completion-irreducible} every irreducible substructure of
$\str{C}$ belongs to $\K^<$.
\end{enumerate}
Then there is $\str{D}\in\K^<$ on the same vertex set as $\str{C}$
such that the identity $\str{C}\to \str{D}$ is a homomorphism-embedding.
\end{proposition}
\begin{proof}
By condition~\ref{cond:completion-order}, choose a linear order
$\prec$ extending $<^{\str{C}}$. Call a pair of distinct vertices
$\{x,y\}$ \emph{ordered} if $x<^{\str{C}}y$ or $y<^{\str{C}}x$.

We claim that the relations of $\str{C}$ specify a partial
alternating function $\omega^\str C$ with values in $\mathbb Z_2$, defined
on precisely those injective $n$-tuples whose vertices
are pairwise ordered by $<^\str{C}$.
Indeed, a tuple in $R^{\str{C}}$ with repeated entries would make
its support an irreducible substructure not belonging to $\K^<$,
contrary to condition~\ref{cond:completion-irreducible}.
If $R^{\str{C}}$ contains a tuple on an $n$-set $X$, then
$\str{C}|_X$ is irreducible, because the tuple contains every
pair of vertices of $X$. Thus $\str{C}|_X\in\K^<$, and all pairs
of $X$ are ordered. Conversely, if all pairs of an $n$-set
$X$ are ordered, then $\str{C}|_X$ is irreducible. By condition~\ref{cond:completion-irreducible},
it belongs to $\K^<$, so its orientation is fully specified by
$R^{\str{C}}$. The claim follows.

We now complete the partial orientation $\omega^\str C$ to
$\omega\colon C^{\underline n}\to\mathbb Z_2$.
Choose a linear order $\sqsubseteq$ on $\binom C2$.
For every $n$-element set $W\subseteq C$ whose orientation
is not yet specified, let $\{u,v\}$, with $u\prec v$, be
the $\sqsubseteq$-least pair of $W$ not ordered in $\str C$.
Enumerate the remaining vertices of $W$ as
$w_1\prec\cdots\prec w_{n-2}$, and put
\begin{equation}\label{eq:finish}
 \omega(u,v,w_1,\ldots,w_{n-2})=0.
\end{equation}
Extend these values by alternation, and put
$\str D=(C,\omega,\prec)$.

Suppose, for a contradiction, that some $\str F\in\F$
embeds into $\str D$, and let $X$ be the image of this
embedding. If every pair of $X$ was already ordered in
$\str C$, then $\str C|_X\in\K^<$, a contradiction.
Otherwise let $\{u,v\}$, with $u\prec v$, be the
$\sqsubseteq$-least pair of $X$ not ordered in $\str C$.
For each $\prec$-increasing $(n-2)$-tuple
$(w_1,\ldots,w_{n-2})$ from $X\setminus\{u,v\}$,
the pair $\{u,v\}$ is also the $\sqsubseteq$-least $<^\str C$-unordered
pair of the set $\{u,v,w_1,\ldots,w_{n-2}\}$.
Thus~\eqref{eq:finish} gives
$\omega(u,v,w_1,\ldots,w_{n-2})=0$.
Consequently, the link of $u,v$ on $X$ is constant, contradicting
pair-cancellation of $\str F$. Hence $\str D\in\K^<$.
Since $\omega^\str C\subseteq \omega$ and $<^\str C\subseteq \prec$, the identity is indeed a homomorphism-embedding $\str C\to \str D$.
\end{proof}

The following proposition was first proved by Sullivan~\cite[Proposition~2.5]{Sullivan2026}. We include a proof for completeness.
\begin{proposition}\label{prop:link-SAP}
Let $\F$ be a (possibly infinite) family of finite
pair-cancelling $n$-hypertournaments. Then both $\K=\Forb(\F)$
and its expansion $\K^<$ by free linear orders have the strong
amalgamation property.
\end{proposition}
\begin{proof}
First consider an amalgamation diagram
$\alpha_i\colon\str{A}\to\str{B}_i$, $i=1,2$, in $\K^<$.
Replacing the structures by isomorphic copies, we may assume
that both maps are inclusions and that
$B_1\cap B_2=A$.
Let $\str{C}$ be the union of $\str{B}_1$ and
$\str{B}_2$. Since the two linear orders agree on $\str{A}$,
$<^{\str{C}}$ admits a linear extension. No relation tuple
of $\str{C}$ contains both a vertex of
$\str{B}_1\setminus\str{A}$ and a vertex of
$\str{B}_2\setminus\str{A}$.
Thus every irreducible substructure of $\str{C}$ lies entirely
in one of the two sides and belongs to $\K^<$.
Proposition~\ref{prop:link-completion} gives the required
strong amalgam.

For unordered inputs, first choose an order on the common
substructure and extend it to both sides; then forget the order
of their completed amalgam.
\end{proof}

\subsection{Ramsey expansions}\label{sec:ramsey-links}
We are now ready to prove Theorem~\ref{thm:intro-link}.
\begin{proof}[Proof of Theorem~\ref{thm:intro-link}]
Let $\mathcal R$ be the class of finite $L$-structures in which
$<$ is a linear order and $R$ is arbitrary. This class is Ramsey
by the Ne\v set\v ril--R\"odl theorem, and all its members are
irreducible.
By Proposition~\ref{prop:link-SAP}, $\K^<$ has strong amalgamation.

To verify local finiteness (Definition~\ref{def:locally-finite}),
fix $\str{C}_0\in\mathcal R$ and let $\str{C}$ satisfy
\ref{item:LF1} and \ref{item:LF2}.
By \ref{item:LF1}, there is a homomorphism-embedding
$\str{C}\to\str{C}_0$, and thus $<^\str C$ has a linear extension.
By \ref{item:LF2}, every irreducible substructure of
$\str{C}$ belongs to $\K^<$.
Consequently, Proposition~\ref{prop:link-completion} gives a strong
completion on the same vertex set, without even using \ref{item:LF3}. We may therefore take $n(\str{C}_0)=1$,
and Theorem~\ref{thm:HN} gives the Ramsey property.

\medskip

The ordering property follows by the standard argument: Fix $\str{A}\in\K$. By the joint embedding property for
$\K^<$, choose $\str{B}^<\in\K^<$ containing every linear
order expansion of $\str{A}$. Since $n\ge3$, there is a unique ordered
$2$-vertex structure $\str{P}^<$, so the Ramsey property gives $\str{C}^<\in\K^<$ with
\[
 \str{C}^<\longrightarrow(\str{B}^<)^{\str{P}^<}_2.
\]
Let $<$ be the order of $\str{C}^<$, and let $\prec$ be any other
linear order of its unordered reduct. Colour each pair according to whether
$<$ and $\prec$ agree or disagree on it. There is a monochromatic
copy of $\str{B}^<$. On this copy, $\prec$ is either $<$ or its reverse.
By our choice of $\str{B}^<$, its order contains all ordered
expansions of $\str{A}$, and so does its reverse.
This proves the ordering property.
\end{proof}

We also get the following corollary:
\begin{corollary}\label{cor:sullivan-ramsey}
For every $n\ge3$ there are $2^{\aleph_0}$ nonisomorphic countably
infinite homogeneous $n$-hypertournaments whose ages have Ramsey
expansions by free linear orders with the ordering property.
\end{corollary}
\begin{proof}
Sullivan~\cite[Proposition~4.11]{Sullivan2026} constructed, for each
$n\ge3$, an infinite embedding-antichain
$\{\str{F}_i:i\in\mathbb N\}$ of finite pair-cancelling
$n$-hypertournaments. For $S\subseteq\mathbb N$, put
$\K_S=\Forb(\{\str{F}_i:i\in S\})$.
By Proposition~\ref{prop:link-SAP}, each $\K_S$ has strong
amalgamation. Since it contains the empty and one-vertex
structures, each $\K_S$ is unbounded, and its
Fra\"iss\'e limit is countably infinite.
The antichain ensures that the classes $\K_S$ are
pairwise distinct, so their Fra\"iss\'e limits are
nonisomorphic. Finally, Theorem~\ref{thm:intro-link}
gives the claimed Ramsey expansions and ordering property.
\end{proof}

\begin{remark}\label{rem:bicycles}
In~\cite[Theorem~2]{cherlin2021ramsey} we constructed arbitrarily
large minimal obstructions to completing unordered partial (\textit{holey})
$\str{H}_4$-free $3$-hypertournaments. We subsequently inferred that the freely ordered class
$\K^<$ could not be a locally finite subclass of $\mathcal R$.
This inference was incorrect. (The original non-completion
construction remains valid.)

The reason is that the old obstructions violate
\ref{item:LF2}. More concretely, consider the first two
implication gadgets in Section~3 of~\cite{cherlin2021ramsey}:
$123\Rightarrow234\Rightarrow345$. The orientation of $123$
is unspecified, while those of $134$, $142$ and $253$
are prescribed. These triples contain the pairs
$13$, $12$ and $23$, respectively.
If the partial structure in the ordered language satisfied
\ref{item:LF2}, each of these three triples, being
irreducible, would itself belong to $\K^<$.
Hence all three pairs of $\{1,2,3\}$ would be compared
by the order relation. The induced substructure on
$\{1,2,3\}$ would then be irreducible, yet its triple
would have no orientation, contradicting~\ref{item:LF2}.
\end{remark}

\section{Even hypertournaments}\label{sec:potential}
Let $\str T = (T,\phi)$ be a tournament (here, $\phi$ is its orientation map). Define a $3$-hypertournament $(T,\omega)$ by $\omega(x,y,z)=\phi(x,y)+\phi(x,z)+\phi(y,z)$.
Observe that, for every $\{x_1,x_2,x_3,x_4\} \in \binom{T}{4}$, it holds that
$\sum_{i=1}^4 \omega(x_1,\ldots,\widehat{x_i},\ldots,x_4) = 0$.
(This holds because every edge of $\str T$ occurs twice in the sum.) It turns out that this extends to every odd arity.

Throughout this section, $n\ge3$ will be odd.
We will use $\omega$ for alternating $n$-ary functions,
$\phi,\psi$ for alternating $(n-1)$-ary functions,
and $s$ for symmetric functions on $(n-2)$-tuples
(or equivalently, functions on $(n-2)$-subsets).

An $n$-hypertournament $\str{A}=(A,\omega)$
is \emph{even} if
\[
 \sum_{i=1}^{n+1}
 \omega(x_1,\ldots,\widehat{x_i},\ldots,x_{n+1})=0
\]
on every injective $(n+1)$-tuple (with the sum
taken in $\mathbb Z_2$). 
We will denote by $\mathcal E_n$ the class of all finite even
$n$-hypertournaments. We will also say that $\omega$ is even.

For a function $f\colon X^{\underline k}\to\mathbb Z_2$,
define 
$\delta f\colon X^{\underline{k+1}}\to\mathbb Z_2$ by
\[
 (\delta f)(x_1,\ldots,x_{k+1})
   =\sum_{i=1}^{k+1}
      f(x_1,\ldots,\widehat{x_i},\ldots,x_{k+1}).
\]
For functions defined on $k$-element subsets,
we use the same notation after identifying them
with symmetric functions on injective tuples.
Note that $\delta^2f=0$, as every term in $\delta^2f$ occurs twice, and that $\omega$ is even if and only if $\delta\omega = 0$.

Suppose that $\phi\colon X^{\underline{n-1}}
\to\mathbb Z_2$ is alternating.
Then $\delta\phi$ is alternating and hence it
is the orientation map of an even $n$-hypertournament.
Indeed, swapping the first two entries gives
\begin{align*}
 &(\delta\phi)(x_2,x_1,x_3,\ldots,x_n)\\
 &\quad=
 \phi(x_1,x_3,\ldots,x_n)
 +\phi(x_2,x_3,\ldots,x_n)\\
 &\qquad+
 \sum_{i=3}^n
 \bigl(\phi(x_1,x_2,x_3,\ldots,\widehat{x_i},
                 \ldots,x_n)+1\bigr)\\
 &\quad=
 (\delta\phi)(x_1,\ldots,x_n)+(n-2)\\
 &\quad=
 (\delta\phi)(x_1,\ldots,x_n)+1.
\end{align*}
The same holds for any transposition of adjacent entries.

It turns out that every even $n$-hypertournament arises in this way, and any two
representing $(n-1)$-hypertournaments differ by switching.
We will now make this precise.

A \emph{potential} of an even $n$-hypertournament
$\str{A}=(A,\omega)$ is an alternating
$(n-1)$-ary function $\phi$ on $A$ satisfying
\[
                    \omega=\delta\phi.
\]

Let $s\colon\binom A{n-2}\to\mathbb Z_2$ be a function.
Note that $\delta s$ is a symmetric function on
$(n-1)$-tuples (and thus $\phi+\delta s$ is
alternating) and $\delta(\phi+\delta s)=\delta\phi$.
We call $\phi\mapsto\phi+\delta s$
\emph{switching} by $s$. Equivalently, switching reverses the
orientation of an $(n-1)$-set $F$ exactly when $F$ contains
an odd number of $(n-2)$-sets $E$ with $s(E)=1$.

\begin{lemma}\label{lem:alternating-sums}
Every finite even $n$-hypertournament has a potential.
Moreover, if $\phi,\psi$ are alternating
$(n-1)$-ary functions on the same vertex set $A$ and
$\delta\phi=\delta\psi$, then there is $s\colon\binom A{n-2}\to\mathbb Z_2$ such that $\psi=\phi+\delta s$.
\end{lemma}
\begin{proof}
Let $\str{A}=(A,\omega)\in\mathcal E_n$.
Choose a linear
order $<$ on $A$ and let $v$ be its least vertex.
For every increasing $(n-1)$-tuple
$x_1<\cdots<x_{n-1}$, put
\[
 \phi(x_1,\ldots,x_{n-1})=
 \begin{cases}
   0,&x_1=v,\\
   \omega(v,x_1,\ldots,x_{n-1}),&x_1\ne v.
 \end{cases}
\]
Extend $\phi$ to non-increasing tuples by alternation. We will now prove that $\delta\phi=\omega$.
If $v=x_1<\cdots<x_n$, then
$(\delta\phi)(x_1,\ldots,x_n)
 =\phi(x_2,\ldots,x_n)=\omega(x_1,\ldots,x_n)$.
If $v<x_1<\cdots<x_n$, evenness of $\omega$ gives
\[
 \begin{aligned}
 0
 &= (\delta\omega)(v,x_1,\ldots,x_n)\\
 &= \omega(x_1,\ldots,x_n)
    +\sum_{i=1}^{n}
       \omega(v,x_1,\ldots,\widehat{x_i},\ldots,x_n)\\
 &= \omega(x_1,\ldots,x_n)
    +(\delta\phi)(x_1,\ldots,x_n).
 \end{aligned}
\]
Hence $\delta\phi=\omega$ on all increasing
$n$-tuples, and thus on all tuples,
since both functions are alternating. In other words, $\phi$ is a potential of $\omega$.

To see the ``moreover'' part, fix $\phi$ and $\psi$. 
Define $d\colon\binom A{n-1}\to\mathbb Z_2$ by
\[
      d(F)=\phi(\bar F)+\psi(\bar F),
\]
where $\bar F$ is an arbitrary enumeration of $F$.
This is well-defined because the sum of two
alternating functions is symmetric. Since
$\delta\phi=\delta\psi$, we have $\delta d=0$.
Pick an arbitrary $v\in A$ and define
$s\colon\binom A{n-2}\to\mathbb Z_2$ by
\[
 s(E)=
 \begin{cases}
   0,&v\in E,\\
   d(E\cup\{v\}),&v\notin E
 \end{cases}.
\]
As above, one checks that $\delta s=d$, and thus $\psi=\phi+\delta s$.
\end{proof}

Let $\mathcal E_n^*$ be the class of finite ordered
relational structures $\str{A}^*=(A,\omega,\phi,<)$, where $\phi$ is an arbitrary alternating
$(n-1)$-ary function on $A$, $\omega=\delta\phi$, and $<$ is an arbitrary
linear order on $A$.

\begin{proposition}\label{prop:potential-SAP}
For every odd $n\ge3$, both $\mathcal E_n^*$
and $\mathcal E_n$ have strong amalgamation.
\end{proposition}
\begin{proof}
First consider an amalgamation diagram
$\str{A}^*\subseteq\str{B}_1^*,\str{B}_2^*$
in $\mathcal E_n^*$, with
\[
 \str{B}_i^*=(B_i,\omega_i,\phi_i,<_i)
 \qquad(i=1,2).
\]
We may assume that $B_1\cap B_2=A$.
The functions $\phi_1,\phi_2$ and the orders
$<_1,<_2$ agree on $A$. Put $C=B_1\cup B_2$,
and choose a linear order $\prec$ on $C$
extending both orders. Extend $\phi_1\cup \phi_2$
arbitrarily to an alternating $(n-1)$-ary
function $\phi$ on $C$.
This is well-defined because $\phi_1$ and $\phi_2$ agree on $A$. Put $\omega=\delta\phi$.
Then $(C,\omega,\phi,\prec)$ is a strong amalgam
in $\mathcal E_n^*$.

Now suppose that
$\str{A}\subseteq\str{B}_1,\str{B}_2$
is an amalgamation diagram in $\mathcal E_n$.
Choose potentials $\phi_1$ and $\phi_2$ of
$\str B_1$ and $\str B_2$, respectively. Since
$\delta\phi_1$ and $\delta\phi_2$ agree on $A$,
Lemma~\ref{lem:alternating-sums} gives us
$s\colon\binom A{n-2}\to\mathbb Z_2$
with $\phi_2|_A=\phi_1|_A+\delta s$.
Extend $s$ to $\binom{B_2}{n-2}$, and replace
$\phi_2$ by $\phi_2+\delta s$ (using the extension
of $s$). The potentials now agree on $A$.
Choose compatible linear orders and apply the
strong amalgamation of $\mathcal E_n^*$.
Forgetting the extra relations gives a strong
amalgam in $\mathcal E_n$.
\end{proof}

\subsection{Ramsey expansions}\label{sec:ramsey-potentials}
Write $\mathcal U_{n-1}$ for the class of all
finite $(n-1)$-hypertournaments, and
$\mathcal U_{n-1}^<$ for its expansion by
free linear orders. The Ramsey property of $\mathcal U_{n-1}^<$ follows from the classical
Ne\v set\v ril--R\"odl theorem~\cite{nevsetvril1977partitions}.
The ordering property follows from Theorem~\ref{thm:intro-link}
with an empty forbidden family for $n\ge5$, and from
\cite[Theorem~1.3]{AlonPachSolymosi2001} for $n=3$.

\begin{proof}[Proof of Theorem~\ref{thm:potential}]
Every $\str{A}^*=(A,\omega,\phi,<)\in\mathcal E_n^*$
is determined by its \emph{potential reduct}
$\str{A}^-=(A,\phi,<)\in\mathcal U_{n-1}^<$.
Conversely, every $(A,\phi,<)\in\mathcal U_{n-1}^<$
has a unique expansion
$(A,\delta\phi,\phi,<)\in\mathcal E_n^*$.
Since these correspondences preserve embeddings,
$\mathcal E_n^*$ is Ramsey.

To prove the expansion property, fix
$\str{A}=(A,\omega)\in\mathcal E_n$.
Choose an ordered $(n-1)$-hypertournament
$\str{D}^<=(D,\phi_D,<^D)$ containing the
potential reducts of all expansions of $\str A$.
Such a $\str D^<$ exists by the joint embedding property of $\mathcal U_{n-1}^<$.

Let $\str{P}^<$ be the unique ordered
$(n-1)$-hypertournament on $n-2$ vertices.
Using the Ramsey property of $\mathcal U_{n-1}^<$,
choose $\str{C}^<=(C,\phi_C,<^C)$ such that
\[
  \str{C}^<\longrightarrow
  (\str{D}^<)^{\str{P}^<}_2.
\]
By the ordering property of $\mathcal U_{n-1}^<$,
choose an $(n-1)$-hypertournament $\str S=(S,\phi_0)$ such that
for every linear order $\prec$ on $S$, the ordered
structure $\str S^\prec=(S,\phi_0,\prec)$ contains
an ordered copy of $\str C^<$.
Put $\str B=(S,\delta\phi_0)\in\mathcal E_n$.

Consider any expansion
$\str B^*=(S,\delta\phi_0,\psi,\prec)\in\mathcal E_n^*$
of $\str B$.
By Lemma~\ref{lem:alternating-sums},
$\psi=\phi_0+\delta s$ for some
$s\colon\binom S{n-2}\to\mathbb Z_2$.
Choose an embedding $f\colon\str C^<\to\str S^\prec$.
Colour each $(n-2)$-subset $E$ of $C$ by $s(f(E))$.
Embeddings of $\str P^<$ into $\str C^<$ correspond
precisely to these subsets, so our Ramsey arrow gives
an embedding $g\colon\str D^<\to\str C^<$ and
$b\in\mathbb Z_2$ such that
$s(f(g(E)))=b$ for every $E\in\binom D{n-2}$.

Put $h=f\circ g$. For every $F\in\binom D{n-1}$,
\[
 (\delta s)(h(F))
 =\sum_{E\in\binom F{n-2}}s(h(E))
 =(n-1)b=0,
\]
since $n-1$ is even. Hence $\psi=\phi_0$ on
$(n-1)$-tuples from $h(D)$. As $h$ preserves
$\phi_0$ and the linear order, it embeds
$\str D^<$ into $(S,\psi,\prec)$.

Every potential reduct of an expansion of $\str A$
embeds into $\str D^<$ and hence into $(S,\psi,\prec)$.
Such an embedding also preserves the original
$n$-ary orientation, since $\omega=\delta\phi$.
Thus every expansion of $\str A$ embeds into
$\str B^*$, as required.
\end{proof}

\section{The extension property for partial automorphisms}
\label{sec:eppa}

For a finite structure $\str{A}$, an \emph{EPPA-witness} is a finite
structure $\str{B}$ containing $\str A$ as a substructure such that every isomorphism
between induced substructures of $\str{A}$ (a \emph{partial automorphism} of $\str A$) extends to an
automorphism of $\str{B}$. A class has the \emph{extension property
for partial automorphisms} (EPPA, also called the
\emph{Hrushovski property}), if each member has such a witness
within the class.

Hrushovski proved EPPA for finite graphs~\cite{hrushovski1992},
a key ingredient in the work of Hodges, Hodkinson, Lascar, and Shelah on the small index property of the random
graph~\cite{hodges1993b}. The area has seen a lot of progress since~\cite{Herwig1995,herwig1998,herwig2000,hodkinson2003,solecki2005,Conant2015,
Hubicka2018metricEPPA,Konecny2018b,eppatwographs,Konecny2019a,Siniora,Evans3,
Hubicka2018EPPA,BradleyEPPAnumbers,EtedadialiabadiGao2025,HubickaSemigeneric,
Aranda2017,BrennerHubickaKonecny2026};
see~\cite{hubicka2025twenty} for a survey.

Interestingly, EPPA is still open for finite tournaments. This problem was already posed by Herwig and Lascar in 2000, who proved that it is equivalent to a statement about the pro-odd topology on free groups~\cite{herwig2000}. 
Historically, non-EPPA examples beyond classes with
order-like obstructions have been relatively scarce.
Huang, Pawliuk, Sabok, and Wise disproved EPPA for certain structures in an infinite relational language which they called $L$-hypertournaments; their proof went through disproving an equivalent free-group profinite-topology criterion~\cite{Sabok}. These are not
the alternating-orientation hypertournaments considered
in this paper. More recently, reducts of the random
graph involving global complementation were shown not
to have EPPA by Brenner, Hubička, and the author~\cite{BrennerHubickaKonecny2026}.

\medskip

In this paper we contribute uncountably many non-EPPA classes which are not order-like. The following proposition is the group-theoretic core behind the results.

\begin{proposition}\label{prop:higharity-eppa-obstruction}
Let $n\ge4$, let $S$ be an $n$-element set, and let
$\str{A}$ be an $n$-hypertournament on
$S\cup\{x,y\}$, where $x,y\notin S$ are distinct.
Enumerate $S\cup\{x\}$ as $(v_0,\ldots,v_n)$
with $v_n=x$, and let $t_0,\ldots,t_n$ be the
values defined in Section~\ref{sec:local-types}
for $\str{A}[S\cup\{x\}]$.
Suppose that
\begin{enumerate}[label=(\roman*),nosep]
\item the map fixing $S$ pointwise and sending
$x$ to $y$ is an isomorphism
$\str{A}[S\cup\{x\}]\to\str{A}[S\cup\{y\}]$;
\item at least four of $t_0,\ldots,t_{n-1}$
are equal (this condition is independent of the enumeration).
\end{enumerate}
Then no finite $n$-hypertournament is an EPPA-witness for $\str{A}$.
\end{proposition}

The proof uses the following standard lemma on
coprime actions of finite groups.

\begin{lemma}\label{lem:klein-centralizer-generation}
Suppose that $\mathbb Z_2^2$ acts by automorphisms
on a finite group $K$ of odd order. Then
$K$ is generated by its elements
fixed by some nonzero element of $\mathbb Z_2^2$.
\end{lemma}
\begin{proof}
This is the standard coprime-action generation theorem
\cite[Theorem~8.3.4(b)]{KurzweilStellmacher2004},
applied to the noncyclic abelian group $\mathbb Z_2^2$
acting on the odd-order group $K$.
\end{proof}

\begin{proof}[Proof of Proposition~\ref{prop:higharity-eppa-obstruction}]
Choose four distinct vertices
$W=\{a,b,c,d\}\subseteq S$ corresponding to
equal values among $t_0,\ldots,t_{n-1}$.
The permutations $(ab)(cd)$ and $(ac)(bd)$,
fixing $S\setminus W$ and $x$, are automorphisms
of $\str A[S\cup\{x\}]$: they are even and
preserve the values $t_i$. Thus
\[
 E=\langle (ab)(cd),(ac)(bd)\rangle
 \leq \operatorname{Aut}(\str A[S\cup\{x\}])
\]
is isomorphic to $\mathbb Z_2^2$, and every
nonidentity element of $E$ acts on $W$
as a double transposition.
By assumption, the map
\[
 p\colon S\cup\{x\}\longrightarrow S\cup\{y\},
 \qquad p|_S=\mathrm{id}_S,\qquad p(x)=y,
\]
is another partial automorphism of $\str{A}$.

Suppose that $\str{B}$ is a finite EPPA-witness
for $\str{A}$, and put
\[
 K=\operatorname{Aut}(\str{B})_{(S)},
\]
the pointwise stabiliser of $S$. The group
$K$ has odd order (an involution in $K$ exchanging vertices $u,v\notin S$ changes the orientation of any $n$-set formed by
$u,v$ and $n-2$ vertices of $S$).

Let $\widetilde G$ be the subgroup of the
setwise stabiliser of $S$ whose restriction
to $S$ belongs to $E$. EPPA implies that
the restriction $\widetilde G_x\to E$ is
onto, since every element of $E$ extends
to an automorphism fixing $x$.

Choose a Sylow $2$-subgroup $\widehat E$
of $\widetilde G_x$. Since the kernel of
the restriction to $E$ is contained in
the odd-order group $K$, the group
$\widehat E$ has four elements and its
restriction to $S$ is isomorphic to $E$.
In particular, $\widehat E$ fixes $x$,
normalises $K$, and its nonidentity elements
act on $W$ as double transpositions.

We claim that every $k\in K$ commuting with
some $1\ne e\in\widehat E$ must fix $x$.
Indeed, write the action of $e$ on $W$ as
$(uv)(wz)$, where $\{u,v,w,z\}=W$, and
suppose that $ke=ek$. Since $e$ fixes $x$,
it also fixes $kx$, because
$e(kx)=k(ex)=kx$. If $kx\ne x$, then the $n$-set
\[
 (S\setminus\{u,v\})\cup\{x,kx\}
\]
is preserved by $e$, which acts on it
as the single transposition $(wz)$.
This would add $1$ to its orientation,
a contradiction.

By Lemma~\ref{lem:klein-centralizer-generation},
applied to the conjugation action of $\widehat E$
on the odd-order group $K$, these elements $k$
generate $K$. Since they all fix $x$, every
element of $K$ fixes $x$. But an automorphism
extending $p$ belongs to $K$ and maps $x$
to $y\ne x$, a contradiction.
\end{proof}

This implies the following theorem:
\begin{theorem}\label{thm:higharity-no-eppa}
Let $n\ge4$, and put $q=\lfloor(n+1)/2\rfloor$.
\begin{enumerate}[label=(\roman*),nosep]
\item If $n\ge6$, no unbounded hereditary strong
amalgamation class of finite $n$-hypertournaments
has EPPA.
\item No hereditary strong amalgamation class of
finite $n$-hypertournaments containing
$\str{T}_{n,q-1}$ has EPPA. More generally, the
same holds if the class contains $\str{T}_{n,d}$
for some $d\le n-3$.
\item If $\F$ is any family of finite pair-cancelling
$n$-hypertournaments, then $\Forb(\F)$ does not
have EPPA. In fact, there is a single $(n+2)$-vertex
$n$-hypertournament belonging to every such class
which has no EPPA-witness even in the unrestricted
class.
\end{enumerate}
\end{theorem}
\begin{proof}
Recall the values $t_i=s_i+i$ defined in
Section~\ref{sec:local-types} for an enumeration
of any $(n+1)$-vertex $n$-hypertournament.
First, suppose that a hereditary strong amalgamation
class $\K$ contains an $(n+1)$-vertex structure
$\str{D}$ with four equal values $t_i$.
Choose the four corresponding vertices and a vertex
$x$ outside them, and put $S=D\setminus\{x\}$.
Strongly amalgamate two copies of $\str{D}$
over $S$, with distinct extra vertices $x,y$.
The resulting $(n+2)$-vertex structure belongs
to $\K$ and satisfies
Proposition~\ref{prop:higharity-eppa-obstruction},
so it has no EPPA-witness even in the unrestricted
class.

For (i), such a $\str{D}$ exists in every
unbounded hereditary class when $n\ge6$,
as at least four of the $n+1$ values $t_i$ agree.
For (ii), take $\str{D}=\str{T}_{n,d}$:
by the definition of $d$, at least
$n+1-d\ge4$ of the values $t_i$ agree.
Since $q-1\le n-3$, this includes
$\str{T}_{n,q-1}$.

For (iii), let $\mathcal P_n$ be the family
of all finite pair-cancelling $n$-hypertournaments,
and put $\K_0=\Forb(\mathcal P_n)$.
By Proposition~\ref{prop:link-SAP}, $\K_0$
has strong amalgamation, and by
Lemma~\ref{lem:local-types} it contains
$\str{T}_{n,q-1}$. Applying (ii) to $\K_0$
gives an $(n+2)$-vertex structure with no
EPPA-witness even in the unrestricted class.
Since $\K_0\subseteq\Forb(\F)$ for every
family $\F$ in the statement, the same
structure works for all of them.
\end{proof}

\medskip

The even hypertournaments also fail EPPA in every odd arity
at least five. In fact, there is a simple reduction to the
unrestricted class one arity lower.

\begin{proposition}\label{prop:even-eppa-transfer}
Let $n\ge3$ be odd. If either $\mathcal E_n$
or the class of all finite $n$-hypertournaments
has EPPA, then the class of all finite
$(n-1)$-hypertournaments has EPPA.
\end{proposition}
\begin{proof}
Let $\str{A}$ be an $(n-1)$-hypertournament and put
$\phi=\omega_{\str{A}}$. Add a new vertex $*$ and define
an $n$-hypertournament $\widehat{\str{A}}$ by
\[
 \omega_{\widehat{\str{A}}}(*,x_1,\ldots,x_{n-1})
     =\phi(x_1,\ldots,x_{n-1})
\]
and, for tuples avoiding $*$, by
\[
 \omega_{\widehat{\str{A}}}(x_1,\ldots,x_n)
     =(\delta\phi)(x_1,\ldots,x_n).
\]
Extend these values to all enumerations by alternation.
The sum of the $n$-face orientations on an
$(n+1)$-tuple containing $*$ is
$\delta\phi+\delta\phi=0$; on a tuple avoiding
$*$ it is $\delta^2\phi=0$.
Thus $\widehat{\str{A}}$ is even. Every partial
automorphism of $\str{A}$ can be extended 
to a partial automorphism of $\widehat{\str{A}}$ fixing $*$.

Under either assumption, there is an EPPA-witness
$\widehat{\str{B}}$ for $\widehat{\str{A}}$.
Put $B=\widehat B\setminus\{*\}$ and define an
$(n-1)$-hypertournament $\str B=(B,\omega_{\str B})$ by
\[
 \omega_{\str B}(b_1,\ldots,b_{n-1})
 =\omega_{\widehat{\str B}}(*,b_1,\ldots,b_{n-1}).
\]
Then $\str B$ contains $\str A$.
For every partial automorphism $p$ of $\str A$,
an automorphism of $\widehat{\str B}$ extending
$p\cup\{*\mapsto *\}$ fixes $*$ and restricts
to an automorphism of $\str B$ extending $p$.
Thus $\str B$ is an EPPA-witness for $\str A$.
\end{proof}

\begin{corollary}\label{cor:even-no-eppa}
For every odd $n\ge5$, the class of finite even
$n$-hypertournaments fails EPPA.
\end{corollary}
\begin{proof}
Apply Proposition~\ref{prop:even-eppa-transfer} and
Theorem~\ref{thm:higharity-no-eppa} in arity $n-1\ge4$,
with $\F=\varnothing$.
\end{proof}

\begin{proof}[Proof of Theorem~\ref{thm:intro-eppa}]
Assertions (i) and (ii) are
Theorem~\ref{thm:higharity-no-eppa}(i) and (iii), respectively.
Assertion (iii) is the case $n=5$ of
Corollary~\ref{cor:even-no-eppa}.
\end{proof}

For $n=3$, Proposition~\ref{prop:even-eppa-transfer}
shows that EPPA for all finite $3$-hypertournaments 
or for all finite even $3$-hypertournaments would imply EPPA for finite
tournaments. Of course, all of these questions remain open.

The class of finite cyclic orders (in the $3$-hypertournament language) fails EPPA.
Indeed, let $a,b,c$ be three cyclically ordered vertices. Since every
two-vertex substructure is trivial, the map $a\mapsto a$, $b\mapsto c$
is a partial automorphism. It cannot extend to an automorphism of any
finite cyclic order containing $a,b,c$, since an automorphism of a finite
cyclic order fixing a vertex is the identity. 

\medskip

To the best of my knowledge, no unbounded amalgamation class of finite
hypertournaments is known to have EPPA. This motivates the following conjecture.
\begin{conjecture}\label{conj:eppa}
For every $n\ge3$, the age of no countably infinite homogeneous
$n$-hypertournament has EPPA.
\end{conjecture}
(Of course, one could also formulate a bolder conjecture allowing $n=2$ as well.)

\section{Concluding remarks and open questions}\label{sec:conclusion}
\subsection{Generalised potential classes}
The potential construction of Section~\ref{sec:potential}
extends to a wider class of labelled structures as follows:

\begin{definition}\label{def:generalised-potential}
Let $n\ge3$, let $Q$ be a finite set with an involution $q\mapsto q^*$, let $G$ be a finite abelian group with an element $a\in G$ (with $a=0$ if $n$ is even), and let $\eta\colon Q\twoheadrightarrow G$ be a surjection satisfying
\[
 \eta(q^*)=a-\eta(q)\qquad(q\in Q).
\]
Let $X$ be a set and $\omega\colon X^{\underline n}\to Q$ a \emph{$Q$-labelling} such that, for every injective $n$-tuple and
every $\sigma\in S_n$,
\[
 \omega(x_{\sigma(1)},\ldots,x_{\sigma(n)})
 =\begin{cases}
   \omega(x_1,\ldots,x_n),&\sigma\text{ is even},\\
   \omega(x_1,\ldots,x_n)^*,&\sigma\text{ is odd}.
  \end{cases}
\]
Define the \emph{effective labelling} $c=\eta\circ\omega$ and assume that
\[
 \sum_{i=1}^{n+1}(-1)^{i-1}
 c(x_1,\ldots,\widehat{x_i},\ldots,x_{n+1})=0
\]
on every injective $(n+1)$-tuple.

We say that $\str X = (X,\omega)$ is a \emph{generalised potential structure}. Embeddings of generalised potential structures are injections which preserve all labels. When considering a
class of these structures, the data $Q,G,a,\eta$
and the involution on $Q$ are fixed.
\end{definition}

We extend the operator $\delta$ from Section~\ref{sec:potential} to functions with values in an arbitrary abelian group $G$ by setting, for every $f\colon X^{\underline{k}} \to G$,
\[
 (\delta f)(x_1,\ldots,x_{k+1})
 =\sum_{i=1}^{k+1}(-1)^{i-1}
  f(x_1,\ldots,\widehat{x_i},\ldots,x_{k+1}).
\]
A \emph{potential} for a generalised potential
$n$-structure $\str X=(X,\omega)$ is a map
$p\colon X^{\underline{n-1}}\to G$ such that $\eta\circ\omega=\delta p$, and
\[
 p(x_{\sigma(1)},\ldots,x_{\sigma(n-1)})
 =\begin{cases}
   p(x_1,\ldots,x_{n-1}),&\sigma\text{ is even},\\
   a-p(x_1,\ldots,x_{n-1}),&\sigma\text{ is odd},
  \end{cases}
\]
for every injective $(n-1)$-tuple
$(x_1,\ldots,x_{n-1})$ and every $\sigma\in S_{n-1}$.

The following theorem generalises Section~\ref{sec:potential}:
\begin{theorem}\label{thm:affine-labelled}
Fix $n,Q,G,a,\eta$ as in
Definition~\ref{def:generalised-potential} and let $\K$ be the class of all such finite generalised potential structures. The following holds:
\begin{enumerate}
    \item Every $\str X\in \K$ has a potential.
    \item The class $\K$ has the strong amalgamation property.
    \item The expansion $\K^*$ of $\K$ obtained by naming a potential
and adding an arbitrary linear order is a precompact Ramsey
expansion with the expansion property.
\end{enumerate}
\end{theorem}
The proof follows essentially the same ideas as
Section~\ref{sec:potential}, but the extra generality requires additional bookkeeping. For this reason, I decided to only cover the hypertournament case in Section~\ref{sec:potential}. A mostly AI-generated supplement giving a readable proof of Theorem~\ref{thm:affine-labelled}, together with a Lean formalization, is available in the supplementary GitHub repository~\cite{Konecny2026GeneralisedPotentialSupplement}.

For $n$ odd, taking $Q=G=\mathbb Z_2$,
$q^*=q+1$, $a=1$, and $\eta$ the identity
recovers the even $n$-hypertournaments and
their potential expansion. The construction also includes
ordinary hypergraphs (fixed labels and the trivial group),
Kay-graphs~\cite{Papadopoulos2024Kay} (symmetric
$\mathbb Z_2$-labels), and multicoloured hypertournaments;
examples of the latter appear in the classification discussed in
Section~\ref{sec:two-coloured-three}.

\subsubsection{Generalised semigeneric structures}
A related but different potential construction arises from
Cherlin's semigeneric tournament~\cite{Cherlin1998}:
Let $G$ be a finite group (not necessarily abelian) and
let $s\colon G\to G$ be an involution satisfying
\[
 s(gh)=s(h)s(1_G)^{-1}s(g).
\]
Choose a finite set $C$ of colours not containing the symbol $E$, with involution
$c\mapsto c^*$ and a surjection
$\phi\colon C\twoheadrightarrow G$ satisfying
$\phi(c^*)=s(\phi(c))$.

A \emph{generalised semigeneric structure} is a structure in the finite binary language $\{E\}\cup\{R_c:c\in C\}$, where $E$ is an equivalence relation and each
ordered pair of vertices $(x,y)$ from distinct $E$-classes has
exactly one colour $c(x,y)\in C$, with
$c(y,x)=c(x,y)^*$. Writing $F=\phi\circ c$,
we require
\[
 F(x_0,y_0)^{-1}F(x_0,y_1)
 =F(x_1,y_0)^{-1}F(x_1,y_1)
\]
whenever $x_0,x_1$ belong to one $E$-class
and $y_0,y_1$ to another.
Allowing $x_0=x_1$ or $y_0=y_1$ adds only tautological instances of the rectangle law.
When $G=\mathbb Z_2$,
$s(g)=g+1$ and $C=G$, this recovers the
semigeneric tournament.

Generalised semigeneric structures also arise from a notion of a potential. However, in this case, the potential is defined between vertices and other equivalence classes.

In an upcoming paper~\cite{Konecny2026SemigenericStructures}, we show that classes of generalised semigeneric structures have strong amalgamation and EPPA.
More interestingly, they witness a strict distinction between
ordinary EPPA, the existence of a dense locally finite
subgroup of the automorphism group of the Fra\"iss\'e
limit, and coherent EPPA. In particular,
semigeneric tournaments have EPPA (this was proved in~\cite{HubickaSemigenericAMUC,HubickaSemigeneric}) but neither of the
stronger properties, answering a question from~\cite{HubickaSemigeneric}.

\subsection{Classification under local constraints}\label{sec:classification-evidence}
Conjecture~\ref{conj:all-homogeneous} asks whether the two constructions
of this paper, together with the cyclic order in arity three, account
for every countably infinite homogeneous hypertournament. Even a weaker
question is open: must every unbounded amalgamation class of finite
$n$-hypertournaments with $n\ge3$ have \emph{strong} amalgamation?
Both constructions in this paper have this property.

The classification becomes more accessible after restricting to $(n+1)$-constrained classes. Recall that
$\str{T}_{n,d}$, $0\le d\le q=\lfloor(n+1)/2\rfloor$, are the types from
Section~\ref{sec:local-types}, and write $\K^{(n)}_S$ for the
$(n+1)$-constrained class allowing precisely the types in $S$.
The following gives the complete answer in this setting, with the positive implications following from this paper.

\begin{theorem}[\cite{Konecny2026ClassificationSupplement}]
\label{thm:constrained-classification}
Let $n\ge4$, put $q=\lfloor(n+1)/2\rfloor$, and let
$S\subseteq\{\str{T}_{n,0},\ldots,\str{T}_{n,q}\}$ contain $\str{T}_{n,q}$.
Then $\K^{(n)}_S$ has the amalgamation property if and only if
\begin{enumerate}[label=(\roman*),nosep]
\item\label{clas:i} $\str{T}_{n,q-1}\in S$, or
\item\label{clas:ii} $n$ is odd and
$S=\{\str{T}_{n,d}:d\equiv q\pmod2\}$.
\end{enumerate}
In each positive case amalgamation is strong.
For $n=2,3$, every unbounded $(n+1)$-constrained amalgamation class
has strong amalgamation.
\end{theorem}

The condition $\str{T}_{n,q}\in S$ is exactly unboundedness: it is necessary
by the finite Ramsey theorem and sufficient because a canonical
linearly oriented hypertournament has only this local type.
Proposition~\ref{prop:link-SAP} gives strong amalgamation for case~\ref{clas:i}, and Proposition~\ref{prop:potential-SAP} for case~\ref{clas:ii}.
In arity 3, there are the four familiar unbounded $4$-constrained classes:
cyclic orders, even $3$-hypertournaments, $\str{H}_4$-free
$3$-hypertournaments, and all $3$-hypertournaments~\cite{cherlin2021ramsey}.

The negative direction of the classification is computer assisted. An AI-generated expository supplement is available in the supplementary GitHub repository~\cite{Konecny2026ClassificationSupplement}, and the repository also contains a Lean formalisation of Theorem~\ref{thm:constrained-classification}. The proof has two parts. It starts with explicit
one-point diagrams witnessing failure of ordinary amalgamation
in small arities. The finite checks are then combined with transfer
lemmas carrying these obstructions from arity $n$ to $n+2$;
two exceptional families require separate arguments, one
explicit and the other probabilistic.

\medskip

At the next local level in arity three, the unbounded
amalgamation classes defined by constraints on induced substructures
of size at most five have also been classified
\cite{Konecny2026ClassificationSupplement}. There are exactly five such classes:
the four just mentioned and $\Forb(\str{P})$, where $\str{P}$ is the five-vertex
$3$-hypertournament on $\{0,1,2,3,4\}$ whose positively oriented
increasing triples are $012,023,034$. They all have strong
amalgamation ($\str{P}$ is pair-cancelling, so $\Forb(\str{P})$ is covered by Theorem~\ref{thm:intro-link}).
The negative direction here is an exhaustive finite verification:
among the $2^{11}=2048$ normalized unbounded choices of allowed
five-vertex types, $37$ explicit one-point non-amalgamation
diagrams rule out all $2043$ negative cases. The validity of
these obstructions, including the failure of ordinary amalgamation,
and their complete coverage are verified in Lean, see~\cite{Konecny2026ClassificationSupplement}.

\medskip

These restricted classifications provide some evidence for
Conjecture~\ref{conj:all-homogeneous}, but, of course, do not settle it.
In fact, increasing the size of the allowed constraints
continues to produce new amalgamation classes.

\begin{observation}\label{obs:strict-constrained-hierarchy}
For every $k>n\ge3$, there is an unbounded strong
amalgamation class of finite $n$-hypertournaments determined
by forbidding a single $k$-vertex $n$-hypertournament.
\end{observation}
\begin{proof}
We construct pair-cancelling $n$-hypertournaments
$\str F_i=(F_i,\omega_i)$ with $F_i=\{1,\ldots,i\}$,
for $i\ge n+1$, such that $\str F_i=\str F_{i+1}[F_i]$.
Start with $\str F_{n+1}\cong\str{T}_{n,0}$, which is
pair-cancelling by Lemma~\ref{lem:local-types}.

Suppose $\str F_i$ is constructed and put
$B=\{1,\ldots,n-1\}$.
Define an $(n-1)$-hypertournament $\str L_i$ on $F_i$
by setting, for every $n\le x\le i$,
\[
 \begin{aligned}
 \omega_{\str L_i}(1,\ldots,n-1)&=n-1,\\
 \omega_{\str L_i}(1,\ldots,\widehat{k},\ldots,n-1,x)
    &=k-1\qquad(1\le k\le n-1),
 \end{aligned}
\]
where the equalities are in $\mathbb Z_2$.
Orient all remaining $(n-1)$-sets arbitrarily.
Observe that $\str L_i[B\cup\{x\}]\cong\str{T}_{n-1,0}$ for every $n\le x\le i$, since all the normalized face values are zero.

Define $\str F_{i+1}$ on $F_{i+1}=F_i\cup\{i+1\}$
by extending $\omega_i$, putting
\[
 \omega_{i+1}(u_1,i+1,u_2,\ldots,u_{n-1})
 =\omega_{\str L_i}(u_1,\ldots,u_{n-1})
\]
for every $(u_1,\ldots,u_{n-1}) \in F_i^{\underline{n-1}}$. 

It remains to check pair-cancellation for pairs
$(v,i+1)$, where $v\in F_i$.
Fix such a $v$ and a linear order $\prec$ on
$F_i\setminus\{v\}$. If $v\in B$, put
$W=B\cup\{n\}$; otherwise put $W=B\cup\{v\}$.
In either case, $v\in W$ and
$\str L_i[W]\cong\str{T}_{n-1,0}$.
Write $W\setminus\{v\}=\{y_1\prec\cdots\prec y_{n-1}\}$.
Since $\str L_i[W]$ has type $\str{T}_{n-1,0}$,
there is $b\in\mathbb Z_2$ such that, for every
$1\le j\le n-1$,
\[
\begin{aligned}
 &\omega_{i+1}(v,i+1,y_1,\ldots,
                       \widehat{y_j},\ldots,y_{n-1})\\
 &\qquad=\omega_{\str L_i}(v,y_1,\ldots,
                       \widehat{y_j},\ldots,y_{n-1})
 =b+j.
\end{aligned}
\]
As $n-1\ge2$, these values are not all equal.
This proves that $\str F_{i+1}$ is pair-cancelling.

For each $k\ge n+1$, Proposition~\ref{prop:link-SAP}
gives strong amalgamation and unboundedness of
$\Forb(\str F_k)$. This class contains every
$n$-hypertournament on fewer than $k$ vertices
but omits $\str F_k$, so it is $k$-constrained
and not $(k-1)$-constrained.
\end{proof}

\subsubsection{Two-coloured $3$-hypertournaments}
\label{sec:two-coloured-three}
A separate computer-assisted classification concerns a variant of 4-constrained
$3$-hypertournaments in which each oriented triple has
one of two colours. These structures can equivalently
be regarded as pairs of ordinary $3$-hypertournaments:
if $o(x,y,z)\in\mathbb Z_2$ records the orientation and
$c(\{x,y,z\})\in\mathbb Z_2$ the colour, the two
orientations are
\[
 R^0=o,\qquad R^1=o+c.
\]
There are $19$ isomorphism types on four vertices.
Among the $393\,216$ unbounded classes defined by
allowing selected four-vertex types, exactly $16\,274$
have amalgamation, and all of these have strong
amalgamation. In fact, $16\,128$ admit a simple block
amalgamation: arrange the two new parts in blocks around
the common part, orient every new mixed triple according
to this block order, and give it one fixed colour. (Note that in the standard hypertournament case, block amalgamation corresponds to the amalgamation procedure implied by Proposition~\ref{prop:link-completion} with $\prec$ putting one new part first, then the common part, and then the other new part.)

The remaining positive classes include several other amalgamation constructions. In particular, they contain
generalised potential classes as in
Definition~\ref{def:generalised-potential}. (Here, $Q$ has size four as there are two orientations and two colours on every triple.)
For example, the class with $Q=G=\mathbb Z_4$, $q^*=1-q$, $a=1$,
and $\eta$ the identity; or the class with $G=\mathbb Z_3$, $a=0$,
and a surjection $\eta$ from $Q$ to
$\mathbb Z_3$, with $\lvert\eta^{-1}[0]\rvert = 2$. The induced involution
on $\mathbb Z_3$ is $g\mapsto-g$; the involution on $Q$ swaps the two
labels over $0$ and swaps the unique labels over $1$ and $-1$.

Curiously, there are also potential-like classes not handled by Theorem~\ref{thm:affine-labelled}. A particularly interesting
example has a switching-poset presentation: start with
a finite poset $\str P=(P,<_P)$ and define a function
$\sigma_P\colon P^{\underline{2}}\to\mathbb Z_3$ by
\[
 \sigma_P(x,y)=
 \begin{cases}
   1,&x<_P y,\\
   -1,&y<_P x,\\
   0,&x\parallel_P y.
 \end{cases}
\]
The effective labels of triples are given by
$\delta\sigma_P$. To obtain the two-coloured hypertournament,
fix an arbitrary linear order on $P$ and, for each increasing triple,
choose an arbitrary raw label in the $\eta$-fibre over its effective
value. Extend these choices to other enumerations using the involution
on $Q$, and then forget the linear order, $\sigma_P$, and the poset
relation. In particular, triples with effective label $0$ have two
possible raw labels.

Requiring the potential to come from a partial
order imposes an additional constraint, so this is a
proper subclass of the corresponding generalised
potential class from Theorem~\ref{thm:affine-labelled}.
In particular, strong amalgamation of this
poset reduct needs a separate argument.

\subsection{Big Ramsey degrees}\label{sec:big-ramsey}
Big Ramsey degrees generalise the Ramsey property, when one wants to find a copy of the
whole Fra\"iss\'e limit. This area has seen a lot of progress
recently; see~\cite{dobrinen2021ramsey,hubicka2024survey}
for background and definitions.
Since all currently known homogeneous hypertournaments admit
Ramsey expansions in finite relational languages, it is
natural to ask about their big Ramsey degrees.

The generic $n$-hypertournament has finite big Ramsey degrees
for every $n\ge2$. Indeed, after fixing an enumeration of its
vertices, it is quantifier-free interdefinable with the
generic $n$-uniform hypergraph equipped with the same
enumeration. Every subcopy of the hypergraph, with the
induced enumeration, again defines a generic
$n$-hypertournament. Since each finite hypertournament
has only finitely many such encodings, finiteness follows
from the corresponding result for generic
hypergraphs~\cite{braunfeld2023big}.
Similarly, by the potential representation used in the
proof of Theorem~\ref{thm:potential}, naming a potential
makes the generic even $n$-hypertournament, for odd $n$,
quantifier-free interdefinable with the generic
$(n-1)$-hypertournament. Finally, the homogeneous cyclic
order is a quantifier-free reduct of $(\mathbb Q,<)$,
so its big Ramsey degrees are finite by Devlin's
theorem~\cite{devlin1979}.

However, once one forbids a nontrivial pair-cancelling
hypertournament, these simple reductions no longer apply.
This leads to the following question.
\begin{question}
Which countably infinite homogeneous hypertournaments have
finite big Ramsey degrees? In particular, which classes
$\Forb(\F)$, for $\F$ a family of pair-cancelling
hypertournaments, have finite big Ramsey degrees?
What about the $\str{H}_4$-free $3$-hypertournaments?
\end{question}

\bibliographystyle{alpha}
\bibliography{ramsey,references_extra}

\newcommand{\etalchar}[1]{$^{#1}$}
\begin{thebibliography}{JLNVTW14}

\bibitem[ABWH{\etalchar{+}}17]{Aranda2017a}
Andres Aranda, David Bradley-Williams, Eng~Keat Hng, Jan Hubi{\v c}ka,
  Miltiadis Karamanlis, Michael Kompatscher, Mat{\v e}j Kone{\v c}n{\'y}, and
  Micheal Pawliuk.
\newblock Completing graphs to metric spaces.
\newblock {\em Electronic Notes in Discrete Mathematics}, 61:53--60, 2017.
\newblock The European Conference on Combinatorics, Graph Theory and
  Applications (EUROCOMB'17).

\bibitem[ABWH{\etalchar{+}}21]{Aranda2017c}
Andres Aranda, David Bradley-Williams, Eng~Keat Hng, Jan Hubi{\v c}ka,
  Miltiadis Karamanlis, Michael Kompatscher, Mat{\v e}j Kone{\v c}n{\'y}, and
  Micheal Pawliuk.
\newblock Completing graphs to metric spaces.
\newblock {\em Contributions to Discrete Mathematics}, 16(2):71--89, 2021.

\bibitem[ABWH{\etalchar{+}}26]{Aranda2017}
Andres Aranda, David Bradley-Williams, Jan Hubi{\v c}ka, Miltiadis Karamanlis,
  Michael Kompatscher, Mat{\v e}j Kone{\v c}n{\'y}, and Micheal Pawliuk.
\newblock Ramsey expansions of metrically homogeneous graphs.
\newblock {\em European Journal of Combinatorics}, 132:104255, 2026.

\bibitem[AH78]{Abramson1978}
Fred~G. Abramson and Leo~A. Harrington.
\newblock Models without indiscernibles.
\newblock {\em Journal of Symbolic Logic}, 43:572--600, 1978.

\bibitem[APS01]{AlonPachSolymosi2001}
Noga Alon, J{\'a}nos Pach, and J{\'o}zsef Solymosi.
\newblock Ramsey-type theorems with forbidden subgraphs.
\newblock {\em Combinatorica}, 21(2):155--170, 2001.

\bibitem[BCdR{\etalchar{+}}24]{braunfeld2023big}
Samuel Braunfeld, David Chodounsk{\'y}, No{\'e} de~Rancourt, Jan Hubi{\v{c}}ka,
  Jamal Kawach, and Mat{\v{e}}j Kone{\v{c}}n{\'y}.
\newblock Big {R}amsey degrees and infinite languages.
\newblock {\em Advances in Combinatorics}, 2024:4, 2024.
\newblock 26pp.

\bibitem[BHK26]{BrennerHubickaKonecny2026}
Sofia Brenner, Jan Hubi{\v{c}}ka, and Mat{\v{e}}j Kone{\v{c}}n{\'{y}}.
\newblock {EPPA} for reducts of the random graph, 2026.
\newblock arXiv:2610.39760.

\bibitem[BNRR18]{bhat2016ramsey}
Vindya Bhat, Jaroslav Nešetřil, Christian Reiher, and Vojtěch Rödl.
\newblock A {R}amsey class for {S}teiner systems.
\newblock {\em Journal of Combinatorial Theory, Series A}, 154:323--349, 2018.

\bibitem[BP10]{Bodirsky2010}
Manuel Bodirsky and Diana Piguet.
\newblock Finite trees are {R}amsey with respect to topological embeddings.
\newblock arXiv:1002.1557, 2010.

\bibitem[Bra26]{Sam}
Samuel Braunfeld.
\newblock Ramsey expansions of {$\Lambda$}-ultrametric spaces.
\newblock {\em European Journal of Combinatorics}, 132:104257, 2026.

\bibitem[BWCHK25]{BradleyEPPAnumbers}
David Bradley-Williams, Peter~J. Cameron, Jan Hubička, and Matěj Konečný.
\newblock {EPPA} numbers of graphs.
\newblock {\em Journal of Combinatorial Theory, Series B}, 170:203--224, 2025.

\bibitem[Che98]{Cherlin1998}
Gregory Cherlin.
\newblock {\em The Classification of Countable Homogeneous Directed Graphs and
  Countable Homogeneous $N$-tournaments}.
\newblock Number 621 in Memoirs of the American Mathematical Society. American
  Mathematical Society, 1998.

\bibitem[Che19]{Cherlin2019Notes}
Gregory Cherlin.
\newblock Unpublished notes, 2019.

\bibitem[Che22]{Cherlin2022}
Gregory Cherlin.
\newblock {\em Homogeneous Ordered Graphs, Metrically Homogeneous Graphs, and
  Beyond, Volume II: 3-Multi-Graphs and 2-Multi-Tournaments}, volume~54 of {\em
  Lecture Notes in Logic}.
\newblock Cambridge University Press, 2022.

\bibitem[CHKN21]{cherlin2021ramsey}
Gregory Cherlin, Jan Hubi{\v{c}}ka, Mat{\v{e}}j Kone{\v{c}}n{\'y}, and Jaroslav
  Ne{\v{s}}et{\v{r}}il.
\newblock Ramsey expansions of 3-hypertournaments.
\newblock In Jaroslav Ne{\v{s}}et{\v{r}}il, Guillem Perarnau, Juanjo Ru{\'e},
  and Oriol Serra, editors, {\em Extended Abstracts EuroComb 2021}, pages
  696--701. Springer International Publishing, 2021.

\bibitem[Con19]{Conant2015}
Gabriel Conant.
\newblock Extending partial isometries of generalized metric spaces.
\newblock {\em Fundamenta Mathematicae}, 244:1--16, 2019.

\bibitem[Dev79]{devlin1979}
Denis Devlin.
\newblock {\em Some partition theorems and ultrafilters on $\omega$}.
\newblock PhD thesis, Dartmouth College, 1979.

\bibitem[Dob23]{dobrinen2021ramsey}
Natasha Dobrinen.
\newblock Ramsey theory of homogeneous structures: current trends and open
  problems.
\newblock In {\em International Congress of Mathematicians}, pages 1462--1486.
  EMS Press, 2023.

\bibitem[EG25]{EtedadialiabadiGao2025}
Mahmood Etedadialiabadi and Su~Gao.
\newblock On extensions of partial isometries.
\newblock In C.~Aiswarya, Prabal~Kumar Sen, and Shashi~Mohan Srivastava,
  editors, {\em Logic and Its Applications}, volume 15402 of {\em Lecture Notes
  in Computer Science}, pages 3--33. Springer, 2025.

\bibitem[EHKN20]{eppatwographs}
David~M. Evans, Jan Hubi{\v{c}}ka, Mat{\v {e}}j Kone{\v {c}}n{\'{y}}, and
  Jaroslav Ne\v{s}et\v{r}il.
\newblock E{P}{P}{A} for two-graphs and antipodal metric spaces.
\newblock {\em Proceedings of the American Mathematical Society},
  148:1901--1915, 2020.

\bibitem[EHN21]{Evans3}
David~M. Evans, Jan Hubi{\v{c}}ka, and Jaroslav Ne{\v{s}}et{\v{r}}il.
\newblock Ramsey properties and extending partial automorphisms for classes of
  finite structures.
\newblock {\em Fundamenta Mathematicae}, 253:121--153, 2021.

\bibitem[Fra53]{Fraisse1953}
Roland Fra{\"\i}ss{\'e}.
\newblock Sur certaines relations qui g\'en\'eralisent l'ordre des nombres
  rationnels.
\newblock {\em Comptes Rendus de l'Academie des Sciences}, 237:540--542, 1953.

\bibitem[Fra86]{Fraisse1986}
Roland Fra{\"\i}ss{\'e}.
\newblock {\em Theory of relations}.
\newblock Studies in logic and the foundations of mathematics. North-Holland,
  1986.

\bibitem[Her95]{Herwig1995}
Bernhard Herwig.
\newblock Extending partial isomorphisms on finite structures.
\newblock {\em Combinatorica}, 15(3):365--371, 1995.

\bibitem[Her98]{herwig1998}
Bernhard Herwig.
\newblock Extending partial isomorphisms for the small index property of many
  $\omega$-categorical structures.
\newblock {\em Israel Journal of Mathematics}, 107(1):93--123, 1998.

\bibitem[HHLS93]{hodges1993b}
Wilfrid Hodges, Ian Hodkinson, Daniel Lascar, and Saharon Shelah.
\newblock The small index property for $\omega$-stable $\omega$-categorical
  structures and for the random graph.
\newblock {\em Journal of the London Mathematical Society}, 2(2):204--218,
  1993.

\bibitem[HJKS19]{HubickaSemigenericAMUC}
Jan Hubi\v{c}ka, Colin Jahel, Mat{\v e}j Kone{\v c}n{\'y}, and Marcin Sabok.
\newblock Extending partial automorphisms of $n$-partite tournaments.
\newblock {\em Acta Mathematica Universitatis Comenianae}, 88(3):807--811,
  2019.
\newblock Extended abstract for Eurocomb 2019.

\bibitem[HJKS25]{HubickaSemigeneric}
Jan Hubi\v{c}ka, Colin Jahel, Mat{\v e}j Kone{\v c}n{\'y}, and Marcin Sabok.
\newblock Extension property for partial automorphisms of the $n$-partite and
  semigeneric tournaments.
\newblock {\em Combinatorica}, 45(5):50, 2025.

\bibitem[HK26]{hubicka2025twenty}
Jan Hubi{\v{c}}ka and Mat{\v{e}}j Kone{\v{c}}n{\'y}.
\newblock Twenty years of {N}e{\v s}et{\v r}il's classification programme of
  {R}amsey classes.
\newblock {\em Computer Science Review}, 59:100814, 2026.

\bibitem[HKN19]{Hubicka2018metricEPPA}
Jan Hubi{\v{c}}ka, Mat{\v{e}}j Kone{\v{c}}n{\'y}, and Jaroslav
  Ne{\v{s}}et{\v{r}}il.
\newblock A combinatorial proof of the extension property for partial
  isometries.
\newblock {\em Commentationes Mathematicae Universitatis Carolinae},
  60(1):39--47, 2019.

\bibitem[HKN21]{Hubicka2017sauerconnant}
Jan Hubi{\v{c}}ka, Mat{\v{e}}j Kone{\v{c}}n{\'y}, and Jaroslav
  Ne{\v{s}}et{\v{r}}il.
\newblock Conant's generalised metric spaces are {R}amsey.
\newblock {\em Contributions to Discrete Mathematics}, 16(2):46--70, 2021.

\bibitem[HKN22]{Hubicka2018EPPA}
Jan Hubi{\v{c}}ka, Mat{\v{e}}j Kone{\v{c}}n{\'{y}}, and Jaroslav
  Ne{\v{s}}et{\v{r}}il.
\newblock All those {E}{P}{P}{A} classes (strengthenings of the
  {H}erwig--{L}ascar theorem).
\newblock {\em Transactions of the American Mathematical Society},
  375(11):7601--7667, 2022.

\bibitem[HKNS26]{HubickaSauerSmetric}
Jan Hubi{\v{c}}ka, Mat{\v{e}}j Kone{\v{c}}n{\'y}, Jaroslav
  Ne{\v{s}}et{\v{r}}il, and Norbert Sauer.
\newblock Ramsey expansions of {$S$}-metric spaces.
\newblock arXiv:2610.26167, 2026.

\bibitem[HL00]{herwig2000}
Bernhard Herwig and Daniel Lascar.
\newblock Extending partial automorphisms and the profinite topology on free
  groups.
\newblock {\em Transactions of the American Mathematical Society},
  352(5):1985--2021, 2000.

\bibitem[HN17]{Hubicka2017designs}
Jan Hubi{\v{c}}ka and Jaroslav Ne\v{s}et\v{r}il.
\newblock Ramsey theorem for designs.
\newblock {\em Electronic Notes in Discrete Mathematics}, 61:623--629, 2017.
\newblock The European Conference on Combinatorics, Graph Theory and
  Applications (EUROCOMB'17).

\bibitem[HN18]{Hubivcka2014}
Jan Hubi{\v{c}}ka and Jaroslav Ne{\v{s}}et{\v{r}}il.
\newblock Bowtie-free graphs have a {R}amsey lift.
\newblock {\em Advances in Applied Mathematics}, 96:286--311, 2018.

\bibitem[HN19]{Hubicka2016}
Jan Hubi{\v{c}}ka and Jaroslav Ne\v{s}et\v{r}il.
\newblock All those {R}amsey classes ({R}amsey classes with closures and
  forbidden homomorphisms).
\newblock {\em Advances in Mathematics}, 356:106791, 2019.

\bibitem[HO03]{hodkinson2003}
Ian Hodkinson and Martin Otto.
\newblock Finite conformal hypergraph covers and {G}aifman cliques in finite
  structures.
\newblock {\em Bulletin of Symbolic Logic}, 9(03):387--405, 2003.

\bibitem[HPSW19]{Sabok}
Jingyin Huang, Micheal Pawliuk, Marcin Sabok, and Daniel~T Wise.
\newblock The {H}rushovski property for hypertournaments and profinite
  topologies.
\newblock {\em Journal of the London Mathematical Society}, 100(3):757--774,
  2019.

\bibitem[Hru92]{hrushovski1992}
Ehud Hrushovski.
\newblock Extending partial isomorphisms of graphs.
\newblock {\em Combinatorica}, 12(4):411--416, 1992.

\bibitem[HZ25]{hubicka2024survey}
Jan Hubi{\v{c}}ka and Andy Zucker.
\newblock A survey on big {R}amsey structures.
\newblock {\em Zbornik Radova (Beograd)}, 22(30):317--355, 2025.
\newblock Selected topics in combinatorial analysis II.

\bibitem[JLNVTW14]{Jasinski2013}
Jakub Jasi{\'n}ski, Claude Laflamme, Lionel Nguyen Van~Th{\'e}, and Robert
  Woodrow.
\newblock Ramsey precompact expansions of homogeneous directed graphs.
\newblock {\em The Electronic Journal of Combinatorics}, 21(4):\#P4.42, 2014.

\bibitem[Kon19]{Konecny2018b}
Mat{\v e}j Kone{\v c}n{\'y}.
\newblock Semigroup-valued metric spaces.
\newblock Master's thesis, Charles University, 2019.
\newblock arXiv:1810.08963.

\bibitem[Kon20]{Konecny2019a}
Mat{\v e}j Kone{\v c}n{\'y}.
\newblock Extending partial isometries of antipodal graphs.
\newblock {\em Discrete Mathematics}, 343(1):111633, 2020.

\bibitem[Kon23]{Konecny2023phd}
Mat{\v e}j Kone{\v c}n{\'y}.
\newblock {\em Model theory and extremal combinatorics}.
\newblock PhD thesis, Charles University, 2023.

\bibitem[Kon26a]{Konecny2026ClassificationSupplement}
Mat{\v{e}}j Kone{\v{c}}n{\'{y}}.
\newblock Classification of {$(n+1)$}-constrained {$n$}-hypertournaments.
\newblock Mostly AI-generated technical supplement to \emph{On homogeneous
  hypertournaments}; accompanying Lean formalization; PDF:
  \url{https://github.com/shulda/hypertournaments/blob/main/papers/hypertournaments_classification_supplement.pdf},
  2026.

\bibitem[Kon26b]{Konecny2026GeneralisedPotentialSupplement}
Mat{\v{e}}j Kone{\v{c}}n{\'{y}}.
\newblock Generalised potential structures.
\newblock Mostly AI-generated proof supplement to \emph{On homogeneous
  hypertournaments}; accompanying Lean formalization; PDF:
  \url{https://github.com/shulda/hypertournaments/blob/main/papers/generalised_potential_structures_proof.pdf},
  2026.

\bibitem[Kon26c]{Konecny2026SemigenericStructures}
Mat{\v{e}}j Kone{\v{c}}n{\'{y}}.
\newblock {Semigeneric Structures: Separating EPPA, Coherent EPPA, and Dense
  Locally Finite Subgroups}.
\newblock Manuscript in preparation, 2026.

\bibitem[KPT05]{Kechris2005}
Alexander~S. Kechris, Vladimir~G. Pestov, and Stevo Todor{\v c}evi{\' c}.
\newblock Fra{\"\i}ss{\'e} limits, {R}amsey theory, and topological dynamics of
  automorphism groups.
\newblock {\em Geometric and Functional Analysis}, 15(1):106--189, 2005.

\bibitem[KS04]{KurzweilStellmacher2004}
Hans Kurzweil and Bernd Stellmacher.
\newblock {\em The Theory of Finite Groups: An Introduction}.
\newblock Universitext. Springer, 2004.

\bibitem[Lac84]{lachlan1984countable}
Alistair~H. Lachlan.
\newblock Countable homogeneous tournaments.
\newblock {\em Transactions of the American Mathematical Society},
  284(2):431--461, 1984.

\bibitem[Lee26]{Lee2026CRelations}
Hyoyoon Lee.
\newblock Ramsey properties of trees and {$C$}-relations via the partite
  construction, 2026.
\newblock arXiv:2610.18402.

\bibitem[LNVTS10]{laflamme2010partition}
Claude Laflamme, Lionel Nguyen Van~Thé, and Norbert~W. Sauer.
\newblock Partition properties of the dense local order and a colored version
  of {M}illiken’s theorem.
\newblock {\em Combinatorica}, 30(1):83--104, 2010.

\bibitem[Ne{\v{s}}07]{Nevsetvril2007}
Jaroslav Ne{\v{s}}et{\v{r}}il.
\newblock Metric spaces are {R}amsey.
\newblock {\em European Journal of Combinatorics}, 28(1):457--468, 2007.

\bibitem[NR75]{nevsetril1975type}
Jaroslav Ne{\v{s}}etril and Vojtech R{\"o}dl.
\newblock Type theory of partition properties of graphs.
\newblock {\em Recent advances in graph theory (Proceedings Second Czechoslovak
  Sympos., Prague, 1974),(Academia, Prague, 1975)}, pages 405--412, 1975.

\bibitem[NR76]{Nevsetvril1976}
Jaroslav Ne{\v{s}}et{\v{r}}il and Vojt{\v{e}}ch R{\"o}dl.
\newblock The {R}amsey property for graphs with forbidden complete subgraphs.
\newblock {\em Journal of Combinatorial Theory, Series B}, 20(3):243--249,
  1976.

\bibitem[NR77a]{nevsetvril1977partitions}
Jaroslav Ne{\v{s}}et{\v{r}}il and Vojt{\v{e}}ch R{\"o}dl.
\newblock Partitions of finite relational and set systems.
\newblock {\em Journal of Combinatorial Theory, Series A}, 22(3):289--312,
  1977.

\bibitem[NR77b]{Nevsetvril1977}
Jaroslav Ne{\v{s}}et{\v{r}}il and Vojt{\v{e}}ch R{\"o}dl.
\newblock A structural generalization of the {R}amsey theorem.
\newblock {\em Bulletin of the American Mathematical Society}, 83(1):127--128,
  1977.

\bibitem[NR78]{Nesetvril1978}
Jaroslav Ne{\v{s}}et{\v{r}}il and Vojt{\v{e}}ch R{\"o}dl.
\newblock On a probabilistic graph-theoretical method.
\newblock {\em Proceedings of the American Mathematical Society},
  72(2):417--421, 1978.

\bibitem[NVT13]{The2013}
Lionel Nguyen Van~Th{\'e}.
\newblock More on the {K}echris--{P}estov--{T}odorcevic correspondence:
  Precompact expansions.
\newblock {\em Fundamenta Mathematicae}, 222:19--47, 2013.

\bibitem[Pap24]{Papadopoulos2024Kay}
Aristomenis-Dionysios Papadopoulos.
\newblock {\em Model-Theoretic Considerations in Finite Combinatorics}.
\newblock PhD thesis, University of Leeds, 2024.
\newblock Section 3.9.

\bibitem[RR26]{reiher2023girth}
Christian Reiher and Vojt{\v{e}}ch R{\"o}dl.
\newblock {\em The girth {R}amsey theorem}.
\newblock Annals of Mathematics Studies. Princeton University Press, Princeton,
  NJ, 2026.
\newblock In press, arXiv:2308.15589.

\bibitem[Sok13]{sokic2013ramsey}
Miodrag Soki{\'c}.
\newblock Ramsey property, ultrametric spaces, finite posets, and universal
  minimal flows.
\newblock {\em Israel Journal of Mathematics}, 194(2):609--640, 2013.

\bibitem[Sok15]{sokic2015semilattices}
Miodrag Soki{\'c}.
\newblock Semilattices and the {R}amsey property.
\newblock {\em The Journal of Symbolic Logic}, 80(4):1236--1259, 2015.

\bibitem[Sol05]{solecki2005}
S{\l}awomir Solecki.
\newblock Extending partial isometries.
\newblock {\em Israel Journal of Mathematics}, 150(1):315--331, 2005.

\bibitem[SS20]{Siniora}
Daoud Siniora and S{\l}awomir Solecki.
\newblock Coherent extension of partial automorphisms, free amalgamation, and
  automorphism groups.
\newblock {\em The Journal of Symbolic Logic}, 85(1):199--223, 2020.

\bibitem[Sul26]{Sullivan2026}
Rob Sullivan.
\newblock Ultrahomogeneous hypertournaments, 2026.
\newblock arXiv:2611.07292.

\end{thebibliography}

\end{document}